\documentclass[11pt]{article}
\usepackage[margin=1in]{geometry}
\usepackage{amsmath,amssymb,amsthm}
\usepackage{booktabs}
\usepackage{multicol}
\usepackage{xcolor}
\usepackage{hyperref}

\newcommand{\AGL}{\mathrm{A\Gamma L}}
\newcommand{\hash}[1]{{\scriptsize\texttt{#1}}}
\title{New lower bounds for binary constant-weight codes from prescribed-automorphism search, and two quantum codes meeting their upper bounds}
\author{Chris Dzoba\\ \small\texttt{chris@dzoba.com}}
\date{15 October 2026}
\begin{document}
\maketitle

\begin{abstract}
We give explicit binary constant-weight codes that raise the lower bounds on $A(n,6,w)$ in eight cells of Brouwer's
table, all now listed there. The first round (9 October 2026) gave four values: $A(27,6,11)\ge 19305$, $A(27,6,12)\ge 24453$, $A(28,6,10)\ge 18642$ and
$A(28,6,9)\ge 11154$, together with an explicit $A(24,6,12)$ code of 5566 words, now in the table's list of lost
codes. The second round (14 October 2026) gave seven values, listed in the table since 15 October 2026: $A(29,6,10)\ge 27183$, $A(29,6,11)\ge 42042$, $A(29,6,9)\ge 14781$,
$A(28,6,10)\ge 19110$, $A(27,6,12)\ge 24804$, $A(28,6,9)\ge 11388$ and $A(27,6,13)\ge 27783$, three of
which (19110, 24804, 11388) replace first-round values. Every code is a union of orbits of the semilinear affine group
$\AGL(1,27)$, its index-2 subgroup, or $\mathrm{AGL}(1,23)$, acting on all but at most two coordinates. The orbits
were chosen by maximum-weight independent set search on the orbit conflict graph: branch and bound in the first round,
iterated local search in the second. Each code was checked by two separately written programs over all pairs of
words and was regenerated from its generators and orbit representatives. We also report two qubit stabilizer codes,
$[[20,11,4]]$ and $[[43,32,4]]$, which raise the lower bounds of 3 listed at codetables.de to the listed upper bound 4.
\end{abstract}

\section{Introduction}
A binary constant-weight code with parameters $(n,d,w)$ is a set of length-$n$ binary words of Hamming weight $w$ at
pairwise Hamming distance at least $d$; $A(n,d,w)$ is the largest size of such a code. Brouwer maintains a table of
the best known bounds \cite{Brouwer}, building on \cite{BSSS90,BBMOS78}; upper bounds for $n\le 28$ there are
from \cite{AVZ00} unless labelled otherwise. Each lower bound carries a label naming its source.

A standard way to find large codes is to prescribe a permutation group $G$ on the coordinates and look only for
$G$-invariant codes \cite{NKO97,BHLO19}. A $G$-invariant code is a union of $G$-orbits of weight-$w$ words, so the
search runs over orbits instead of words. This note reports codes of this kind for $d=6$ and $27\le n\le 29$ that raise
the lower bound in eight cells: four round-1 values and seven round-2 values, three of which replace round-1 values. All of
them are now listed in the table (checked 15 October 2026), where the three replaced round-1 values no longer appear. It also reports one explicit code for a lost-code entry. All are lower bounds from explicit codes; none
is claimed to be optimal, and we claim improvement only relative to the table entries as checked, not priority over
unpublished work. Bounded literature checks made on 9 October for the round-1 cells found no competing construction in recent heuristic work on the table \cite{Ros26,Lys26,Ech26}; such checks are not proofs of absence.

\section{Results}\label{sec:results}

\paragraph{Groups.} Coordinates $0,\dots,26$ are identified with $\mathrm{GF}(27)=\mathrm{GF}(3)[x]/(x^3+2x+1)$, the element
$c_0+c_1x+c_2x^2$ being coordinate $c_0+3c_1+9c_2$. $\AGL(1,27)$ is the group of maps $z\mapsto az^{3^s}+b$ with
$a\ne 0$, $s\in\{0,1,2\}$; it has order 2106. Its index-2 subgroup $H$ (order 1053) requires $a$ to be a nonzero
square. $\mathrm{AGL}(1,23)$ acts on $\mathrm{GF}(23)$ (coordinates $0,\dots,22$) by $z\mapsto az+b$, order 506. Extra coordinates
(27, 28; or 23 for $n=24$) are fixed. Table~\ref{tab:groups} lists the six group files; orders and point orbits
were recomputed for this note by full closure of the generators.

\begin{table}[h]
\centering\small
\begin{tabular}{llrll}
\toprule
name & group & order & point orbits & file \\
\midrule
$G_{24}$ & $\mathrm{AGL}(1,23)$, $\infty$ fixed & 506 & $23+1$ & \texttt{group-AGL1-23-plus-fixed.txt}\\
$G_{27}$ & $\AGL(1,27)$ & 2106 & $27$ & \texttt{group-AGammaL1-27-on-27.txt}\\
$G_{28}$ & $\AGL(1,27)$, 1 fixed & 2106 & $27+1$ & \texttt{group-AGammaL1-27-plus-fixed.txt}\\
$H_{28}$ & $H$, 1 fixed & 1053 & $27+1$ & \texttt{group-AGammaL1-27-index2-plus-fixed.txt}\\
$G_{29}$ & $\AGL(1,27)$, 2 fixed & 2106 & $27+1+1$ & \texttt{group-AGammaL1-27-plus-2-fixed.txt}\\
$H_{29}$ & $H$, 2 fixed & 1053 & $27+1+1$ & \texttt{group-AGammaL1-27-index2-plus-2-fixed.txt}\\
\bottomrule
\end{tabular}
\caption{The prescribed groups. Generators are listed in Appendix~A.1.}
\label{tab:groups}
\end{table}

\paragraph{Round 1 (codes sent to the table maintainer on 9 October 2026).} The previous values are those of the live
table as checked on 9 October 2026 between 04:06 and 06:54 UTC. In the table, \emph{Bo} is Boisclair (emails of
22 and 30 August 2026), \emph{g} denotes a code with a prescribed group, and \emph{BJL\"O} is \cite{BHLO19}. The
18260 entry carried no label; it was marked as first given on the table page.

\begin{table}[h]
\centering\small
\begin{tabular}{lrrrlrrl}
\toprule
cell & previous [label] & new & gain & group & orbits & $V$ & search \\
\midrule
$A(27,6,11)$ & 18260 [--] & 19305 & $+1045$ & $G_{27}$ & 13 & 2227 & B\&B completed, 95.1 s \\
$A(27,6,12)$ & 23903 [Bo] & 24453 & $+550$ & $G_{27}$ & 18 & 3648 & B\&B cap 1500 s \\
$A(28,6,10)$ & 18346 [Bo] & 18642 & $+296$ & $G_{28}$ & 14 & 1217 & B\&B completed, 124.2 s \\
$A(28,6,9)$ & 10920 [g] & 11154 & $+234$ & $H_{28}$ & 18 & 1918 & B\&B cap 900 s \\
\midrule
$A(24,6,12)$ & 5558 [BJL\"O]$^{*}$ & 5566 & $+8^{*}$ & $G_{24}$ & 12 & 1788 & B\&B completed, 754.7 s \\
\bottomrule
\end{tabular}
\caption{Round-1 codes. $V$ is the number of vertices of the orbit conflict graph. ${}^{*}$For $A(24,6,12)$ the
table's lower bound remains the 5616 of \cite{BSSS90}, whose code is lost; 5558 was the best explicitly available
code, and 5566 replaces it in the table's list of lost codes. It is not an improved lower bound.}
\label{tab:r1}
\end{table}

For the three completed runs the branch and bound finished without a time cap being reached, so 19305, 18642 and
5566 are the largest codes invariant under $G_{27}$, $G_{28}$ and $G_{24}$ respectively, provided the orbit graphs
are correct. These completions were reported by the solver and were not repeated independently; for the 5566 run the
saved output, graph and witness were checked, but no launch manifest records which binary ran (a later rerun with the
repaired solver reached the same 12 orbits in 580.5 s). The 24453 and 11154 runs stopped at their time caps after
$1{,}233{,}305{,}600$ and $629{,}970{,}944$ search nodes and are lower bounds only.

\paragraph{Round 2 (codes sent to the table maintainer on 14 October 2026).} The previous values are those of the
live table at 08:15 EDT on 14 October 2026, which agree with the first fetch of Section~\ref{sec:status}. \emph{Dz} is
round 1 above; \emph{Gi} is Gippini (email of 1 October 2026, \cite{Gip26}); \emph{gs} is the Graham--Sloane
construction with an $S_2$ set \cite{GS80}; \emph{JX} is \cite{JX14}.

\begin{table}[h]
\centering\footnotesize
\begin{tabular}{lrrrlrrl}
\toprule
cell & previous [label] & new & gain & group & orbits & $V$ & first run reaching it \\
\midrule
$A(29,6,10)$ & 25896 [Gi] & 27183 & $+1287$ & $H_{29}$ & 39 & 7174 & ILS2, seed 2, 900 s \\
$A(29,6,11)$ & 41216 [gs] & 42042 & $+826$ & $G_{29}$ & 28 & 4335 & ILS2, seed 3, 900 s \\
$A(29,6,9)$ & 14313 [Gi] & 14781 & $+468$ & $H_{29}$ & 19 & 2310 & ILS2, seed 1, 240 s \\
$A(28,6,10)$ & 18642 [Dz] & 19110 & $+468$ & $H_{28}$ & 28 & 5256 & ILS2 seeded by 18993, seed 5, 900 s \\
$A(27,6,12)$ & 24453 [Dz] & 24804 & $+351$ & $G_{27}$ & 15 & 3648 & ILS2, seed 3, 1200 s \\
$A(28,6,9)$ & 11154 [Dz] & 11388 & $+234$ & $H_{28}$ & 18 & 1918 & ILS, seed 1, 600 s \\
$A(27,6,13)$ & 27600 [JX] & 27783 & $+183$ & $G_{27}$ & 17 & 4332 & ILS2, seed 1, 240 s \\
\bottomrule
\end{tabular}
\caption{Round-2 codes. ILS and ILS2 are the iterated local searches of Section~\ref{sec:method}; ILS2 used
mean acceptance slack 300 in these runs. Three entries (19110, 24804, 11388) supersede round-1 values:
$19110>18642$, $24804>24453$, $11388>11154$. All seven values are timed lower bounds, not group optima.}
\label{tab:r2}
\end{table}

Later runs reproduced the same values: unseeded ILS2 with other seeds returned 19110, 24804 and 11388 again, and
ILS2 restarted from each final code (one or two chains per cell, 600 to 1200 s each) did not exceed any of the seven. Intermediate
codes, verified the same way and superseded, were $A(29,6,10)\ge 26247, 26598$; $A(29,6,11)\ge 41418, 41691$;
$A(29,6,9)\ge 14430$; and $A(28,6,10)\ge 18993$.

Orbit sizes, from regeneration: 27183 has orbits $39+117+17\cdot 351+20\cdot1053$; 42042 has
$39+234+2\cdot351+6\cdot702+1053+17\cdot2106$; 14781 has $39+6\cdot351+12\cdot1053$; 19110 has
$39+117+12\cdot351+14\cdot1053$; 24804 has $234+2\cdot702+2\cdot1053+10\cdot2106$; 11388 has
$39+117+8\cdot351+8\cdot1053$; 27783 has $54+2\cdot351+702+1053+12\cdot2106$.

\paragraph{SHA-256 of the word lists.} Each word list is a text file of 0/1 rows, character $i$ = coordinate $i$,
one row per word, terminated by newlines (file names \texttt{A$n$-6-$w$-$M$.txt}).

\begin{table}[h]
\centering
\begin{tabular}{ll}
\toprule
code & SHA-256 of word list \\
\midrule
$A(24,6,12)\ge 5566$  & \hash{a330246135644f9d1a1a8f2bb46ffa5865eeafa3b1c64bfa03e073495cea4d07}\\
$A(27,6,11)\ge 19305$ & \hash{0adb0037a0bcb033912eaba346527315b1f4d11f745403c42b2af3d326586943}\\
$A(27,6,12)\ge 24453$ & \hash{2e2c3e4befc59546b19ba7dcf2693f324bbd2d577e6ed9d646a883b92d91bf57}\\
$A(28,6,10)\ge 18642$ & \hash{b23744d89191690ac01eb9d4d122dc3648755e37c0eef4982e45453c027bd28d}\\
$A(28,6,9)\ge 11154$  & \hash{1a3b1dbbaa5b17a099284ce26b81a41f703f1cb22146416391f8448ebf606a02}\\
\midrule
$A(29,6,10)\ge 27183$ & \hash{436d58bbe6fbf03711a79033c0d3b70210fff53affc74bdb8915e32dd223db39}\\
$A(29,6,11)\ge 42042$ & \hash{30a50d6aea7beef477579f79fb9ec09ebd5b8c8dc515dafdcf791b6144347c83}\\
$A(29,6,9)\ge 14781$  & \hash{c0cf7c71f59a524b1618a3e2f4dbae7808057ce309698c9141db26fb41a011e5}\\
$A(28,6,10)\ge 19110$ & \hash{dcc386f01f9fdcfdb9e4779300f279798a0254f43c56ac100d06744bea8f0c03}\\
$A(27,6,12)\ge 24804$ & \hash{817d070244da7316c66b8845b03cb3ea9d84069984b8d4e5619e657c790b9b3c}\\
$A(28,6,9)\ge 11388$  & \hash{7635db5fa873c47a63af70e77e1b076c0fcd39b5bb44c205d2a33dc64017ea1a}\\
$A(27,6,13)\ge 27783$ & \hash{550d7f4f896950d126da3c4b2dfd29a34e8762cec3a919f3b3cbcd950fd1502a}\\
\bottomrule
\end{tabular}
\caption{Hashes recomputed for this note by \texttt{check\_np.py}; they equal the hashes in the packages sent on
9 and 14 October.}
\label{tab:sha}
\end{table}

\section{Table status}\label{sec:status}
The live table \cite{Brouwer} was fetched twice; the author keeps both snapshots.
{\sloppy\begin{itemize}
\item 15 October 2026, 01:50:51 UTC (21:50:51 EDT, 14 October), 127{,}211 bytes, server \texttt{Last-Modified}
14 October 2026 13:36:17 GMT. The round-1 values were listed with label Dz; the seven round-2 target cells still showed
their previous values. The round-2 codes had been emailed to the maintainer at 21:43 EDT on 14 October.
SHA-256:\newline{\tiny\texttt{7e2afb60f7237aeb82a1f1e1f46358341408e4ab7b38f3272da201f8aee4b977}}
\item 15 October 2026, 10:44:13 UTC (06:44:13 EDT), 127{,}299 bytes, server \texttt{Last-Modified} 15 October 2026
07:27:15 GMT. Its key reads ``Dz: Chris Dzoba, email 2026-10-09,
2026-10-15 with 5+7 new codes.'' SHA-256:\newline{\tiny\texttt{7c8870a49f4db267162a08d81f35679e95c1f6172dcf0baf3173805dd9e27cfd}}
\end{itemize}}

\begin{table}[h]
\centering\small
\begin{tabular}{lrll}
\toprule
cell & our value & live table at 06:44 EDT, 15 Oct & status \\
\midrule
$A(27,6,11)$ & 19305 & $19305^{\rm Dz}$--$66078^{\rm KKT}$ & listed \\
$A(28,6,10)$ & 18642 & $19110^{\rm Dz}$--63756 & was listed; replaced by 19110 \\
$A(27,6,12)$ & 24453 & $24804^{\rm Dz}$--$84573^{\rm KKT}$ & was listed; replaced by 24804 \\
$A(28,6,9)$ & 11154 & $11388^{\rm Dz}$--31195 & was listed; replaced by 11388 \\
$A(24,6,12)$ & 5566 & lost codes: 5616 [BSSS], $5566^{\rm Dz}$ & listed (as explicit code) \\
\midrule
$A(29,6,10)$ & 27183 & $27183^{\rm Dz}$--\dots & listed \\
$A(29,6,11)$ & 42042 & $42042^{\rm Dz}$--\dots & listed \\
$A(29,6,9)$ & 14781 & $14781^{\rm Dz}$--\dots & listed \\
$A(28,6,10)$ & 19110 & $19110^{\rm Dz}$--63756 & listed \\
$A(27,6,12)$ & 24804 & $24804^{\rm Dz}$--$84573^{\rm KKT}$ & listed \\
$A(28,6,9)$ & 11388 & $11388^{\rm Dz}$--31195 & listed \\
$A(27,6,13)$ & 27783 & $27783^{\rm Dz}$--$91079^{\rm KKT}$ & listed \\
\bottomrule
\end{tabular}
\caption{Table status at the second fetch (15 October 2026, 06:44 EDT). KKT is \cite{KKT11}; ``\dots'' means no upper
bound is listed. In the main $d=6$ table $A(24,6,12)$ still reads 5616--15906, the 5616 marked as lost.}
\label{tab:status}
\end{table}

\section{Quantum codes}\label{sec:quantum}
The two codes below are qubit stabilizer codes \cite{CRSS98}. The construction and search details come from the
notes of the session that found them. The four stabilizer files in the release match the SHA-256 prefixes recorded there,
and the referee re-checked them for this note with a fresh script (\texttt{qcheck.py}; our rerun is in
\texttt{reruns/qcheck.txt}): for $[[20,11]]$, rank 9, all generators commuting, no logicals of weight $\le 3$ and
1{,}075 of weight 4; for $[[43,32]]$, rank 11, all generators commuting, no logicals of weight $\le 3$, the weight-4
logical $X_0X_1Z_2Y_{14}$, minimum stabilizer weight 22, and $g_1,g_2 \mid y^{21}-1$. At codetables.de
\cite{Grassl}, the cells $[[20,11]]$ and $[[43,32]]$ list lower bound 3 and upper bound 4. Both codes therefore meet
the listed upper bound; we do not claim more than that.

\paragraph{Construction.} Both codes are additive and invariant under a qubit permutation $\sigma$; the stabilizer is
spanned by the $\sigma$-orbits of one or two seed Paulis. Qubits $0,\dots,f-1$ are fixed by $\sigma$; the remaining
$n-f$ qubits form two blocks of length $m$, and $\sigma$ shifts position $p$ to $p+1 \bmod m$ in both blocks. Block
coordinates are read as polynomials in $\mathrm{GF}(2)[y]/(y^m-1)$, separately for the $X$ and $Z$ parts. For each
seed the $X$ part is chosen first and the $Z$ part is any solution of the linear system expressing commutation with
the rows already present and with the seed's own $\sigma^s$-images. To obtain rank exactly $n-k$, the seed's block
polynomials are forced into the ideal generated by a divisor $g$ of $y^m-1$, so that the seed's orbit has dimension
$m-\deg g$.
\begin{itemize}
\item $[[20,11,4]]$: $f=0$, $m=10$, one seed, $g=y+1$; $\sigma$ has two 10-cycles. Nine generators. The code is pure,
with minimum stabilizer weight 12.
\item $[[43,32,4]]$: $f=1$, $m=21$, two seeds, $g_1=1+y^2+y^7+y^8+y^{11}$ and $g_2=(y^{21}-1)/(y-1)$, both divisors of
$y^{21}-1$ over $\mathrm{GF}(2)$; orbit dimensions $10+1=11=n-k$; $\sigma$ has one fixed point and two 21-cycles. The
code is pure, with minimum stabilizer weight 22.
\end{itemize}
They were found by random sampling within this family (programs \texttt{cycsample3} and \texttt{cycsample5}, random
seed 1), at sample 376 for $[[20,11,4]]$ (11 October) and sample 4378 for $[[43,32,4]]$ (12 October), in sweeps over
table entries whose bounds differ by one. An overnight search on 13--14 October over 73 such entries with $n\le 60$
and $d\le 6$ found no further improvement.

\paragraph{Verification.} Each code was checked to have minimum distance exactly 4 by
(A) a Brouwer--Zimmermann enumeration on the length-$3n$ binary image that tests every low-weight normalizer word for
membership in $S$, validated against brute force on small and random codes; (B) a second Brouwer--Zimmermann
implementation written by a separate agent that had not seen (A), with brute-force and meet-in-the-middle modes,
validated on known codes and 200 random codes; and an exhaustive scan of every Pauli of weight 1 to 3 (32{,}550 for
$n=20$, 341{,}463 for $n=43$). The meet-in-the-middle mode of (B) confirmed $d=4$ for $[[20,11,4]]$; it was not run on
$[[43,32,4]]$, which exceeds its limit $n\le 30$. A referee using fresh code found, for $[[20,11,4]]$: rank 9, all
generators commuting, no logical of weight $\le 3$ and 1{,}075 of weight 4, matching the MacWilliams transform of the
stabilizer weight distribution $\{0{:}1, 12{:}50, 14{:}220, 16{:}195, 18{:}20, 20{:}26\}$, and $d=3$ on the table's own
$[[20,11,3]]$ as a negative control. For $[[43,32,4]]$: rank 11, normalizer dimension 75, no logical of weight
$\le 3$, and a full weight-4 scan of 9{,}996{,}210 Paulis giving 5{,}740 logicals, none in $S$ (witness
$X_0X_1Z_2Y_{14}$), consistent with MacWilliams; negative controls gave $d=3$ for $[[5,1,3]]$, $d=2$ with a planted
weight-3 logical, and $d=1$ with a generator dropped. The referee's scripts are in \texttt{referee/q20verify/} and
\texttt{referee/q43verify/} of the ancillary files.

\paragraph{Derivability and status.} At $n=20$, $[[20,12]]$ and $[[20,13]]$ have upper bound 3, and the standard
propagation rules do not take $[[21,11,4]]$ or $[[22,12,4]]$ to $[[20,11,4]]$. At $n=43$, a $[[44,33,4]]$ is listed;
the referee showed that the stored $[[44,33,4]]$ and the stored $[[64,53,4]]$ do not shorten to $[[43,32,4]]$. There
remains a risk that an inequivalent $[[44,33,4]]$, not stored in the table, shortens to a $[[43,32,4]]$. No prior
$[[20,11,4]]$ or $[[43,32,4]]$ was found in the literature searches recorded by that session. Both codes were emailed
to codes@codetables.de, on 11 and 12 October 2026. No reply had been received by 21:55 EDT on 14 October, and a
fetch of both cells at 01:56:21 UTC on 15 October (21:56 EDT, 14 October) shows lower bound 3 and upper bound 4 for
each. A further fetch at 03:07 UTC on 16 October 2026 still showed lower
bound 3 and upper bound 4 for both cells.

\section{Method}\label{sec:method}
\paragraph{Orbit conflict graph.} Fix $G$ and $(n,d,w)$. Compute the $G$-orbits on weight-$w$ words, discard every
orbit containing two words at distance $<d$, and join two remaining orbits by an edge if some word of one is at
distance $<d$ from some word of the other. Each vertex is weighted by its orbit size. $G$-invariant codes with minimum
distance $\ge d$ are exactly the independent sets, and the code size is the total weight. The builder \texttt{og}
stores the orbit of every word in a hash table; the round-2 builder \texttt{og3} replaces it with an array indexed by
colex rank, which needs $4\binom{n}{w}$ bytes (about 310~MB for $\binom{29}{14}$) and produced a graph byte-identical
to \texttt{og}'s for $(29,6,9)$ under the twisted group below.

\paragraph{Exact search.} \texttt{wmis} is a bitset branch and bound for maximum-weight independent sets with a
weighted colouring bound. A completed run gives the maximum over $G$-invariant codes, provided the graph is correct;
a run stopped by its time cap gives a lower bound. Round 1 used it throughout.

\paragraph{Iterated local search.} In round 2 two heuristics were run on the same graphs, giving lower bounds only.
\texttt{ils} uses $(w,1)$-swaps (insert a vertex and drop its neighbours in the solution if the weight increases) with
forced random insertions as perturbation. \texttt{ils2} adds $(1,2)$-swaps, accepts a new local optimum if its weight
is at least the current weight minus a slack drawn from an exponential distribution (mean 300 in the runs of Table~\ref{tab:r2}; later confirmation runs also used 150 and 800),
and periodically restarts from the best solution; \texttt{ils2s} starts it from a verified code. ILS found larger codes than
branch and bound on the same graphs. On $A(29,6,10)$ under $G_{29}$ ($V=1579$), branch and bound reached 26247 in
3000 s, and ILS reached 26598 after 7.6 s. On $A(29,6,11)$ under the twisted group described below ($V=6085$),
branch and bound reached 37830 in 2400 s and ILS2 reached 41418; ILS2 under $G_{29}$ later reached 42042. In a rerun
for this note (\texttt{reruns/ils-repro.txt}, ILS seed 1), ILS reached two exact group optima: 58149
($\mathrm{PSL}(2,27)$ on 28 points plus one fixed, $w=12$) in under 0.1 s and 14430 (the twisted group, $w=9$) in
2.6 s. The random seed mattered: on $A(27,6,12)$ under $G_{27}$, ILS seed 1 (600 s) and ILS2 seed 1 (150 s) both
stopped at 24219, while ILS2 seed 3 reached 24804.

\paragraph{Why the $n=29$ cells had not been searched.} The round-1 job lists were generated from table cells that have an upper bound
(\texttt{tablejobs.py}, \texttt{tablejobs2.py}: \texttt{if not ub: continue}). The cells $A(29,6,8\ldots14)$ and
$A(29\ldots32,8,8)$ list no upper bound, so round 1 never created jobs for them; they were filtered out, not skipped
for memory. The round-1 sweeps of $n=29$ to 32 covered only cells that do list an upper bound. In round 2 the first
120-second branch-and-bound run on $A(29,6,10)$ under $G_{29}$ already exceeded 25896.

\paragraph{The twisted action and the index-2 subgroup.} $\AGL(1,27)$ also acts on $27+2$ points with the quadratic
character of the multiplier swapping the two extra points, while translations and the Frobenius map fix them. Full
closure gives order 2106 and point orbits $27+2$; this group is not a subgroup of $G_{29}$ (file
\texttt{n29\_gfagaml\_q27\_k26tw.txt}). Its exact optimum at $(29,6,9)$ is 14430, above the 14313 then listed, while
the optimum under $G_{29}$ is exactly 14313. The common subgroup of the plain and twisted actions is $H$ with both
extra points fixed, i.e.\ $H_{29}$ (squares and Frobenius, order 1053); every code invariant under either group is $H_{29}$-invariant, and with ILS2 it
gave the largest codes at $A(29,6,9)$, $A(29,6,10)$, $A(28,6,10)$ and $A(28,6,9)$. The same subgroup, with one fixed
point, had given $A(28,6,9)\ge 11154$ in round 1, where the optimum under $G_{28}$ is exactly 10920 (a completed run on
a 362-vertex graph), equal to the value then listed.

\section{Verification}\label{sec:verify}
\paragraph{Two checkers per code.} \texttt{verify.c} (C) reads the 0/1 rows and compares every pair character by
character, reporting rows, wrong weights, minimum distance and violations. \texttt{check\_np.py} (Python/NumPy,
sharing no code with the builders, solvers or \texttt{verify.c}) checks lengths, alphabet, weights and distinctness,
computes every pairwise distance as $2(w-\langle x,y\rangle)$ in blocks, counts the pairs, and recounts the closest
pair by XOR. Both were re-run for this note on all twelve files in the release folder; both report VALID with minimum
distance exactly 6 for every code, and \texttt{check\_np.py} checked exactly $M(M-1)/2$ pairs in each case
(Table~\ref{tab:verify}). Every word list was also rebuilt from its generator file and orbit-representative file
(\texttt{expand\_fast.py}) and equals the stored list as a set, with the orbit counts of Tables~\ref{tab:r1}
and~\ref{tab:r2}.

\begin{table}[h]
\centering\small
\begin{tabular}{lrrl|lrrl}
\toprule
code & $M$ & pairs checked & & code & $M$ & pairs checked & \\
\midrule
$A(24,6,12)$ & 5566 & 15{,}487{,}395 & & $A(29,6,10)$ & 27183 & 369{,}444{,}153 & \\
$A(27,6,11)$ & 19305 & 186{,}331{,}860 & & $A(29,6,11)$ & 42042 & 883{,}743{,}861 & \\
$A(27,6,12)$ & 24453 & 298{,}962{,}378 & & $A(29,6,9)$ & 14781 & 109{,}231{,}590 & \\
$A(28,6,10)$ & 18642 & 173{,}752{,}761 & & $A(28,6,10)$ & 19110 & 182{,}586{,}495 & \\
$A(28,6,9)$ & 11154 & 62{,}200{,}281 & & $A(27,6,12)$ & 24804 & 307{,}606{,}806 & \\
 & & & & $A(28,6,9)$ & 11388 & 64{,}837{,}578 & \\
 & & & & $A(27,6,13)$ & 27783 & 385{,}933{,}653 & \\
\bottomrule
\end{tabular}
\caption{Pair counts from \texttt{check\_np.py} (round 1 left, round 2 right). Every code: 0 wrong weights, minimum distance 6 and 0 violations from both checkers; all rows distinct
(checked by \texttt{check\_np.py}, implied by minimum distance 6 in \texttt{verify.c}).}
\label{tab:verify}
\end{table}

\paragraph{At the time of discovery.} In round 1 each code was expanded from its generators (\texttt{expand.py}) and
checked by \texttt{verify.c}; Codex then rebuilt each group with its own field arithmetic (all 2106, 1053 or 506
maps, compared with a separate generator closure), regenerated the orbits, and compared every pair character by
character, with rejection controls (duplicate, too-close, wrong-weight, wrong-length and wrong-character rows). In
round 2 each code was checked by \texttt{verify.c} and \texttt{check\_np.py}, the latter validated on the round-1 11154
code (62{,}200{,}281 pairs, the same count as Codex's) and on a too-close negative control.

\paragraph{Referee.} A separate Claude session refereed every round-2 code (including the superseded intermediate
ones) with code written afresh: a C checker (\texttt{chk.c}) over all pairs with a full distance histogram; three
negative controls per code (a second row moved to distance 2 from the first, a duplicated row, a row of wrong weight),
each of which must fail; orbit regeneration from the \texttt{.reps} and \texttt{.gens} files by full group closure,
with a check of group order, point orbits, orbit sizes and disjointness; a re-fetch of the live table for each cell;
and checks that no listed entry yields the new values by shortening or extension. It reported every code confirmed
(\texttt{referee/cw/referee.json}).

\paragraph{Limits of these checks.} The codes of Jin and Xing for $A(28,6,12)\ge 40237$ and $A(28,6,13)\ge 49510$ were
not reconstructed, so their shortenings were not computed; a shortening would have to beat the parent's average by 7.9\%
to reach 24804 at $A(27,6,12)$ and by 4.8\% to reach 27783 at $A(27,6,13)$. For $A(29,6,11)$ we did not search the Graham--Sloane
family \cite{GS80}, which uses a Sidon set in an abelian group of order at least $813=29\cdot 28+1$; since the
best class of such a partition can exceed its average, we cannot rule out that some choice beats 42042. Shortenings of the new $n=29$ codes give at most 10569 at $A(28,6,9)$ and at most 18603
at $A(28,6,10)$, below the new values there. The branch-and-bound optima depend on the correctness of the graph
builders, which were cross-checked (\texttt{og} against \texttt{og3} on one graph) but not proved correct.

\section{Negative results}\label{sec:neg}
Round 2 recorded 1043 jobs. Of these, 848 ran: 666 completed exact, 173 timed lower bounds,
8 vertex-limit skips and 1 deliberately stopped run, with no memory skips. The other 195 were deferred by throttle
windows and never run; 82 (cell, group) pairs among them were not rerun. Table~\ref{tab:neg} gives the best value for each main target that did not improve, separating group
maxima (completed branch and bound; valid for that group only, provided the graph is correct) from timed searches.

\begin{table}[h]
\centering\small
\begin{tabular}{lrll}
\toprule
cell & listed & best completed exact (group) & best timed (method) \\
\midrule
$A(27,6,8)$ & 4828 & 4050 (twisted $\mathrm{AGL}(1,25)$ on $25+2$) & -- \\
$A(27,6,9)$ & 8324 & 8112 ($G_{27}$ and $H$, round 1) & 8112 ($H$, ILS2) \\
$A(27,6,10)$ & 13026$^{\dagger}$ & -- & 12987 ($H$, ILS2) \\
$A(27,6,11)$ & 19305 & 19305 ($G_{27}$, round 1) & 19305 ($H$, ILS2 seeded); 18603 unseeded \\
$A(28,6,8)$ & 6356 & 5105 (twisted $\mathrm{AGL}(1,25)$, 1 fixed) & -- \\
$A(28,6,11)$ & 29484 & 24948 ($\mathrm{PSU}(3,3)$ on 28) & 28782 ($G_{28}$, ILS2, 3 seeds) \\
$A(28,6,12)$ & 40237 & 37107 ($\mathrm{PSU}(3,3)$ on 28) & 38259 ($G_{28}$, ILS2) \\
$A(28,6,13)$ & 49510 & 43344 ($\mathrm{PSU}(3,3)$ on 28) & 44982 ($G_{28}$, ILS2) \\
$A(28,6,14)$ & 53046 & 46224 ($\mathrm{PSU}(3,3)$ on 28) & 47439 ($G_{28}$, ILS2) \\
$A(29,6,8)$ & 8148 & 7488 ($\mathrm{AGL}(1,27)$, 2 fixed) & -- \\
$A(29,6,12)$ & 61820 & 58149 ($\mathrm{PSL}(2,27)$, 1 fixed) & 57447 ($H_{29}$, ILS2) \\
$A(29,6,13)$ & 80837 & 70434 ($\mathrm{PSL}(2,27)$, 1 fixed) & 68607 ($\mathrm{PSU}(3,3)$, 1 fixed, 60 s) \\
$A(29,6,14)$ & 92366 & 71892 ($\mathrm{PSL}(2,27)$, 1 fixed) & 77076 ($\mathrm{PSU}(3,3)$, 1 fixed, 60 s) \\
$A(29,8,8)$ & 1060 & $\le 505$ (all tested groups) & -- \\
\bottomrule
\end{tabular}
\caption{Main targets not improved. Exact values are group maxima from completed runs (caps 60--120 s in round 2);
timed values are lower bounds. ${}^{\dagger}$13026 [Fr] since 13 October; 12987 at the round-2 baseline.}
\label{tab:neg}
\end{table}

Other round-2 families stayed well below the table and are all completed-exact negatives unless noted: $W(E_6)$ on
the 27 lines and its rotation subgroup, with fixed or twisted extensions; $\mathrm{Sp}(6,2)$ on the 28 bitangents;
$\mathrm{PSU}(3,3)$ and $\mathrm{P\Gamma U}(3,3)$ on the 28-point Hermitian unital; $\mathrm{PSL}(2,8)$ and
$\mathrm{P\Gamma L}(2,8)$ on 28 points; twisted affine groups over $\mathrm{GF}(23)$ and $\mathrm{GF}(25)$ at $n=25$--29 and
$d=6,8,10$ (60-second triage, some deferred and not rerun); $\mathrm{AGL}(1,29)$ and its order-351 subgroup at $n=29$
(very large graphs, weak timed values).

Round 1 swept the 186 open cells with $n\le 28$ (3105 jobs: 2462 completed exact, 385 timed, 256 skipped, 2 target
hits) and the cells with $29\le n\le 32$ that list an upper bound (798 completed exact, 779 skipped by memory or vertex
gates, 31 timed, no hits), over affine and semilinear groups over prime and prime-power fields with fixed points,
$\mathrm{PSL}/\mathrm{PGL}/\mathrm{P\Gamma L}(2,q)$, products and diagonal actions. These counts are solver-output
classifications, not independently certified exclusions, and early runs with undersized solver arrays are not usable
as negatives. Specific round-1 group maxima: $A(24,6,11)$, 5106 under $G_{24}$ (listed 5376); $A(26,6,9)$, 5600 under an
affine group over $\mathrm{GF}(25)$ with one fixed point (listed 5761); $A(28,8,6)$, 81 for codes invariant under a 27-cycle
fixing one point (listed 132; two graph reconstructions agree). Timed only: $A(27,6,13)$, 26730 under $G_{27}$ in 900 s;
$A(27,6,10)$, 12636 under $H$ in 900 s; $A(23,6,10)$, 2783 under the order-253 subgroup of $\mathrm{AGL}(1,23)$ in 900 s
(listed 2992).

\section{Data availability}
The word lists, orbit representatives, group generators, both checkers, the regeneration script, the checker outputs,
the quantum stabilizers, the referee scripts and a \texttt{SHA256SUMS} file are provided as ancillary files with
this arXiv submission (\texttt{codes/}, \texttt{groups/}, \texttt{checkers/}, \texttt{verification/},
\texttt{quantum/}, \texttt{referee/}, \texttt{reruns/}). Generators and orbit representatives are also printed in
Appendix~A, which suffices to rebuild every code. The table snapshots cited above are not redistributed.

\section*{Tool and computational resource disclosure}
This work was carried out by two AI systems working under the author's direction, and this note discloses their roles
in full. In round 1 (night of 8--9 October 2026) Claude (Anthropic) wrote the orbit-graph builder and the weighted
branch and bound, generated the group families, ran the table sweeps and the longer pushes that produced the five
codes, and checked each hit with its own expansion and checker. Codex (OpenAI) pinned the live table and parsed it
into job targets, traced the provenance of the lost codes, audited the solver (finding a stack-depth defect that was
then repaired, and a too-small array that invalidated some early negatives), reconstructed every hit independently
with its own field arithmetic and a literal pairwise check, computed shortenings, and reviewed the submission package.
In round 2 (night of 13--14 October) a Claude session wrote the array-indexed graph builder, the larger solver, the two
iterated local searches, the twisted-group generator and the second checker, and ran all searches; a separate Claude
session refereed every new code with freshly written programs. The quantum codes were found by a Claude session; their
second distance checker was written by a separate agent that had not seen the first, and they were refereed by
another Claude session. This paper was written by Claude from the project's notes, logs and verification outputs.
The author set the goals, approved the runs, decided what to submit, sent the submissions, and is responsible for the
content. Computing: one laptop (Apple M4 Max, 128 GB); round 1 shared at most four worker processes between the two
agents, and round 2 ran at most three single-threaded processes at a time, all at reduced priority. Software: C and
Python 3 with NumPy; C++ for the quantum search and checkers; \texttt{tectonic} for this document.

\appendix
\section*{Appendix A. Generators, orbit representatives and stabilizers}
\subsection*{A.1 Group generators}
{\scriptsize
\noindent\textbf{G24} (\texttt{group-AGL1-23-plus-fixed.txt}): $\mathrm{AGL}(1,23)$, coordinate 23 fixed, order 506.\\
\begin{verbatim}
1 2 3 4 5 6 7 8 9 10 11 12 13 14 15 16 17 18 19 20 21 22 0 23
0 5 10 15 20 2 7 12 17 22 4 9 14 19 1 6 11 16 21 3 8 13 18 23
\end{verbatim}
\noindent\textbf{G27} (\texttt{group-AGammaL1-27-on-27.txt}): $\mathrm{A\Gamma L}(1,27)$ on 0..26, order 2106.\\
\begin{verbatim}
1 2 0 4 5 3 7 8 6 10 11 9 13 14 12 16 17 15 19 20 18 22 23 21 25 26 24
3 4 5 6 7 8 0 1 2 12 13 14 15 16 17 9 10 11 21 22 23 24 25 26 18 19 20
9 10 11 12 13 14 15 16 17 18 19 20 21 22 23 24 25 26 0 1 2 3 4 5 6 7 8
0 3 6 9 12 15 18 21 24 5 8 2 14 17 11 23 26 20 7 1 4 16 10 13 25 19 22
0 1 2 5 3 4 7 8 6 13 14 12 15 16 17 11 9 10 26 24 25 19 20 18 21 22 23
\end{verbatim}
\noindent\textbf{G28} (\texttt{group-AGammaL1-27-plus-fixed.txt}): $\mathrm{A\Gamma L}(1,27)$, coordinate 27 fixed, order 2106.\\
\begin{verbatim}
1 2 0 4 5 3 7 8 6 10 11 9 13 14 12 16 17 15 19 20 18 22 23 21 25 26 24 27
3 4 5 6 7 8 0 1 2 12 13 14 15 16 17 9 10 11 21 22 23 24 25 26 18 19 20 27
9 10 11 12 13 14 15 16 17 18 19 20 21 22 23 24 25 26 0 1 2 3 4 5 6 7 8 27
0 3 6 9 12 15 18 21 24 5 8 2 14 17 11 23 26 20 7 1 4 16 10 13 25 19 22 27
0 1 2 5 3 4 7 8 6 13 14 12 15 16 17 11 9 10 26 24 25 19 20 18 21 22 23 27
\end{verbatim}
\noindent\textbf{H28} (\texttt{group-AGammaL1-27-index2-plus-fixed.txt}): index-2 subgroup, coordinate 27 fixed, order 1053.\\
\begin{verbatim}
1 2 0 4 5 3 7 8 6 10 11 9 13 14 12 16 17 15 19 20 18 22 23 21 25 26 24 27
3 4 5 6 7 8 0 1 2 12 13 14 15 16 17 9 10 11 21 22 23 24 25 26 18 19 20 27
9 10 11 12 13 14 15 16 17 18 19 20 21 22 23 24 25 26 0 1 2 3 4 5 6 7 8 27
0 9 18 5 14 23 7 16 25 15 24 6 11 20 2 13 22 4 21 3 12 26 8 17 19 1 10 27
0 1 2 5 3 4 7 8 6 13 14 12 15 16 17 11 9 10 26 24 25 19 20 18 21 22 23 27
\end{verbatim}
\noindent\textbf{G29} (\texttt{group-AGammaL1-27-plus-2-fixed.txt}): $\mathrm{A\Gamma L}(1,27)$, coordinates 27, 28 fixed, order 2106.\\
\begin{verbatim}
1 2 0 4 5 3 7 8 6 10 11 9 13 14 12 16 17 15 19 20 18 22 23 21 25 26 24 27 28
3 4 5 6 7 8 0 1 2 12 13 14 15 16 17 9 10 11 21 22 23 24 25 26 18 19 20 27 28
9 10 11 12 13 14 15 16 17 18 19 20 21 22 23 24 25 26 0 1 2 3 4 5 6 7 8 27 28
0 3 6 9 12 15 18 21 24 5 8 2 14 17 11 23 26 20 7 1 4 16 10 13 25 19 22 27 28
0 1 2 5 3 4 7 8 6 13 14 12 15 16 17 11 9 10 26 24 25 19 20 18 21 22 23 27 28
\end{verbatim}
\noindent\textbf{H29} (\texttt{group-AGammaL1-27-index2-plus-2-fixed.txt}): index-2 subgroup, coordinates 27, 28 fixed, order 1053.\\
\begin{verbatim}
1 2 0 4 5 3 7 8 6 10 11 9 13 14 12 16 17 15 19 20 18 22 23 21 25 26 24 27 28
3 4 5 6 7 8 0 1 2 12 13 14 15 16 17 9 10 11 21 22 23 24 25 26 18 19 20 27 28
9 10 11 12 13 14 15 16 17 18 19 20 21 22 23 24 25 26 0 1 2 3 4 5 6 7 8 27 28
0 9 18 5 14 23 7 16 25 15 24 6 11 20 2 13 22 4 21 3 12 26 8 17 19 1 10 27 28
0 1 2 5 3 4 7 8 6 13 14 12 15 16 17 11 9 10 26 24 25 19 20 18 21 22 23 27 28
\end{verbatim}
}
\subsection*{A.2 Orbit representatives}
Each entry is \texttt{mask:size}, with bit $i$ of the mask the value of coordinate $i$.
{\scriptsize
\noindent\textbf{$A(24,6,12)\ge 5566$}, group G24, 12 orbits:
\begin{verbatim}
163065:253 8494451:253 93115:506 8429007:506 8447663:506 304631:506 186015:506 8533467:506
180107:506 8433559:506 8426487:506 73691:506
\end{verbatim}
\noindent\textbf{$A(27,6,11)\ge 19305$}, group G27, 13 orbits:
\begin{verbatim}
828029:702 276441:702 16115:702 338663:702 2367379:702 357345:1053 308202:2106 325843:2106
315189:2106 309622:2106 288506:2106 861172:2106 305343:2106
\end{verbatim}
\noindent\textbf{$A(27,6,12)\ge 24453$}, group G27, 18 orbits:
\begin{verbatim}
4095:234 875406:351 2472277:351 834262:351 113885:351 114099:351 866021:702 2432915:702 842343:2106
837500:2106 863732:2106 351730:2106 844509:2106 382649:2106 316405:2106 310899:2106 292787:2106
812007:2106
\end{verbatim}
\noindent\textbf{$A(28,6,10)\ge 18642$}, group G28, 14 orbits:
\begin{verbatim}
134218239:39 134507910:351 804047:702 308124:702 39417:702 29547:702 134255837:702 273393:2106
274150:2106 313579:2106 277619:2106 134485718:2106 134485817:2106 134231924:2106
\end{verbatim}
\noindent\textbf{$A(28,6,9)\ge 11154$}, group H28, 18 orbits:
\begin{verbatim}
511:39 112420:117 37814:117 39416:351 113740:351 134486484:351 134330885:351 134219486:351
134233868:351 134257091:351 40711:1053 283223:1053 286353:1053 281073:1053 287604:1053 265845:1053
134483426:1053 134256356:1053
\end{verbatim}
\noindent\textbf{$A(29,6,10)\ge 27183$}, group H29, 39 orbits:
\begin{verbatim}
29621:351 48310:1053 65161:351 113601:351 267251:1053 269199:1053 276970:1053 286628:1053
289466:1053 293573:1053 351779:351 373644:351 391720:1053 829583:1053 134218239:39 134223535:1053
134256703:351 134330148:117 134331468:351 134481658:1053 134487657:1053 134495043:351
134502571:1053 134563780:351 268443215:1053 268467793:351 268473565:351 268473710:351 268473779:351
268549410:351 268699622:1053 268708651:1053 268716515:1053 268717166:1053 268727573:351
402655157:351 402667203:351 402692834:1053 402918742:1053
\end{verbatim}
\noindent\textbf{$A(29,6,11)\ge 42042$}, group G29, 28 orbits:
\begin{verbatim}
30571:702 48535:2106 272318:2106 285658:2106 292975:2106 304087:2106 312190:2106 805358:2106
830434:702 2427851:702 134257148:702 134482807:2106 134493052:2106 134510937:2106 134531017:2106
268476199:1053 268710646:2106 268711271:2106 268713874:2106 268715479:2106 268729621:702
402653695:39 402667350:2106 402690998:234 402691438:702 402925755:2106 402943366:351 402970826:351
\end{verbatim}
\noindent\textbf{$A(29,6,9)\ge 14781$}, group H29, 19 orbits:
\begin{verbatim}
511:39 13942:1053 40235:1053 113600:351 266029:1053 271226:1053 282561:1053 290182:351
134222819:1053 134233681:351 134258202:1053 134485829:1053 268437214:351 268483160:1053
268548611:351 268707129:1053 268715625:1053 402658005:1053 402692548:351
\end{verbatim}
\noindent\textbf{$A(28,6,10)\ge 19110$}, group H28, 28 orbits:
\begin{verbatim}
29547:351 39417:351 40906:1053 65298:351 113210:351 113747:1053 269235:1053 271739:1053 275417:351
275956:351 285870:1053 289242:351 290630:1053 305779:1053 337836:1053 448145:351 134218239:39
134221555:351 134250065:351 134256342:1053 134256703:351 134330148:117 134481767:1053
134485406:1053 134486844:1053 134494085:1053 134495043:351 134524556:1053
\end{verbatim}
\noindent\textbf{$A(27,6,12)\ge 24804$}, group G27, 15 orbits:
\begin{verbatim}
4095:234 48827:1053 277855:2106 307966:2106 314350:2106 316349:2106 321462:2106 326547:2106
339897:2106 354261:2106 792503:2106 818629:1053 839079:2106 1873117:702 1873262:702
\end{verbatim}
\noindent\textbf{$A(28,6,9)\ge 11388$}, group H28, 18 orbits:
\begin{verbatim}
511:39 38323:351 39416:351 40711:1053 112420:117 113740:351 265845:1053 281073:1053 283223:1053
286353:1053 287604:1053 134219486:351 134233868:351 134256356:1053 134257091:351 134330885:351
134483426:1053 134486484:351
\end{verbatim}
\noindent\textbf{$A(27,6,13)\ge 27783$}, group G27, 17 orbits:
\begin{verbatim}
360093:2106 375479:2106 386989:1053 448221:351 817567:2106 825599:2106 833387:2106 847566:2106
873342:2106 884133:2106 916262:2106 1874922:2106 1895651:351 2407397:2106 2408348:702 2450679:2106
3974005:54
\end{verbatim}
}
\subsection*{A.3 Quantum stabilizers}
Pauli strings; character $i$ acts on qubit $i$. Release files \texttt{quantum/q*\_paulis.txt}.
{\scriptsize
\noindent\textbf{[[20,11,4]], 9 generators}:
\begin{verbatim}
YZIIZZIZIXZXZXZYZXXX
IYZIIIIZIZYIZYYZZZZI
ZIXZZZZZIXXIZYYIZIXI
ZIIYZIIIIZIZYIZYYZZZ
ZZZIYIZIZYZIYXIIXXYX
ZZZIIYZIIIIZIZYIZYYZ
IIIZIZYIZYYZZZZIIYZI
IIZZIZIXZXZXZYZXXXYZ
ZIIIIZIZYIZYYZZZZIIY
\end{verbatim}
\noindent\textbf{[[43,32,4]], 11 generators}:
\begin{verbatim}
YZIIIIIIIIIXYXYZZZXIXXZYIYYIIZZIYXXIYXXYZYI
IXIZIIIIIIIXYIIYXIYZIXYZZIYYYZZZXXZXZXZZYIZ
ZZXIZZZZZIIZXXXZZXYIYIXIZIIXIYIZYZYXZIZZZII
IIZXIZIIIIIIIXYIIYXIYZIXYZZIYYYZZZXXZXZXZZY
ZIIZXIZZZZZIIZXXXZZXYIYIXIZIIXIYIZYZYXZIZZZ
ZIZIIYZIZIIZXXXYXXYYZIZZZXYYZYYXIIXYYXIYIXZ
ZZZIIZXIZZZZZIIZXXXZZXYIYIXIZIIXIYIZYZYXZIZ
ZIZZZIIYZIZXYZXZIYZYYZXZIZXZIXIXZIZXYZXYYYZ
ZIZZZIIZXIZZZZZIIZXXXZZXYIYIXIZIIXIYIZYZYXZ
ZZZZZZZIIYZXXXYXIZYXZZIIYZYXZIZYXIIXIIXIZYX
ZIZZZZZIIZXXXZZXYIYIXIZIIXIYIZYZYXZIZZZIIZX
\end{verbatim}
}

\end{document}